\documentclass{asl}
\usepackage[section]{placeins}
 
\usepackage{pstricks}
\usepackage{epsf}
\usepackage{amsfonts,amssymb}

\usepackage{graphicx}

\newtheorem{Theorem}{Theorem} 
 
\newtheorem{Lemma}{Lemma}

\psset{unit=3cm}
\def\axioms{\begin{tabular*}{0.8\textwidth}{@{\extracolsep{\fill}}lr}}
\def\endaxioms{\end{tabular*}\smallskip}
\def\implies{\ \rightarrow\ }

\def\R{{\mathbb R}}

\def\ilc{{i\ell c}}
\def\K{{\mathbb K}}
\def\F{{\mathbb F}}
\def\L{{\mathbb L}}
\def\O{{\mathcal O}}

\def\B{{\bf B}}
\def\T{{\bf T}}
\def\E{{\bf E}}

\newcommand{\bT}[3]{\T(\MakeLowercase{#1},\MakeLowercase{#2},\MakeLowercase{#3})}
\newcommand{\congT}[2]{\E(\MakeLowercase{#1},\MakeLowercase{#2})}
\def\imp{\implies}

\def\SegmentCircleContinuityFigure{%
\pspicture(0,0.2)(3,2.2)
\psdot(1.45,1.05)
\put(1.45,0.92){$a$}
\pscircle(1.45,1.05){1.05}
\qline(0,1.20)(2.4,1.8)
\psdot(0.2,1.25)
\put(0.13,1.33){$q$}
\psdot(0.4,1.05)
\put(0.31,0.9){$b$}
\qline(0,1.05)(1.45,1.05)
\psdot(0.65,1.05)
\put(0.7,0.92){$x$}
\psdot(0.18,1.05)
\put(0.08,0.92){$y$}
\psarcn[linestyle=dashed](1.45,1.05){0.80}{190}{125}
\psarcn[linestyle=dashed](1.45,1.05){1.27}{200}{145}
\psdot(0.72, 1.38)
\put(0.74,1.25){$p$}
\pscircle(0.436,1.31){0.035}
\put(0.35,1.35){$z$} 
\endpspicture}

\def\CounterExampleFigure{%
\psset{unit=14mm}
\pspicture(-4,1.1)(4,0)
\psgrid[gridcolor=lightgray,subgriddiv=5,griddots=5, gridlabels=7pt](-4,0)(4,1)
\psdot(0,0)
\put(-0.15,0.04){$t$}
\psline(0,0)(0,1)
\psline(-1,0.75)(1,0.75)
\psdot(0,1)
\put(-0.01,1.07){$a$}
\psdot(0,0.75)
\put(0.05,0.55){$d$}
\psdot(-1,0.75)
\put(-1.18,0.8){$b$}
\psdot(1,0.75)
\put(1.1,0.8){$c$}
\pscircle(-4,0){0.05}
\psline(-4,0)(0,1)
\pscircle(4,0){0.05}
\psline(4,0)(0,1)
\endpspicture}

\def\InnerOuterPaschFigure{%
\pspicture(3.6,1.8)
\psdot(0,0.6)
\put(0,0.47){$a$}
\pscircle(1,0.6){0.03}
\put(1,0.47){$x$}
\qline(0,0.6)(0.97,0.6)
\psdot(1.3,0.6)
\put(1.35,0.55){$q$}
\qline(1.03,0.6)(1.3,0.6)
\psdot(1.1,1.4)
\put(1.1,1.45){$c$}
\qline(0,0.6)(1.1,1.4)  
\qline(1.1,1.4)(1.375,0.3) 
\psdot(1.375,0.3)
\put(1.43,0.25){$b$}
\psdot(0.52,0.98)
\put(0.43,1.03){$p$}
\qline(1.375,0.3)(1.015,0.584)  
\qline(0.52,0.98)(0.979,0.62)  
\psdot(2,0.3)
\put(1.95,0.19){$a$}
\psdot(3.5,0.3)
\put(3.45,0.19){$q$}
\qline(2,0.3)(2.47,0.3)
\qline(2.53,0.3)(3.5,0.3)
\pscircle(2.5,0.3){0.03}
\put(2.45,0.19){$x$}
\psdot(2.7,1.5)
\put(2.65,1.57){$b$}
\qline(2.7,1.5)(3.5,0.3)  
\psdot(3.0,1.05)
\put(3.07,1.02) {$c$}
\qline(2.7,1.5)(2.5,0.33)  
\qline(2,0.3)(3.0,1.05)   
\psdot(2.566,0.722)
\put(2.6,0.6){$p$}
\endpspicture}

\def\TarskiParallelFigure{%
\pspicture(2.5,1.2)
\qline(0.03,0)(2.47,0)  
\pscircle(0,0){0.03}
\put(-0.13,-0.02){$x$}
\pscircle(2.5,0){0.03}
\put(2.55,-0.02){$y$}
\psdot(0.5,1)
\put(0.47,1.05){$a$}
\psdot(0.167,0.333)
\put(0.06,0.33){$b$}
\qline(0.5,1)(2.48,0.022)  
\qline(0.01,0.025)(0.5,1)  
\psdot(1.1,0.703)
\put(1.1,0.74){$c$}
\qline(0.167,0.333)(1.1,0.703)  
\qline(0.5,1)(1.0,0)  
\psdot(1.0,0)
\put(1.02,0.05){$t$}
\psdot(0.723,0.553)
\put(0.78,0.48){$d$}
\endpspicture
}

\def\EuclidParallelFigure{%
\pspicture (0, 0.2)(3.2, 1.2)
\qline(0.0,0.15)(2.47,0.15)
\qline(2.53,0.15)(2.7,0.15)
\qline(0.0,0.7)(2.7,0.7)
\qline(0.0,0.9)(2.473,0.163)
\qline(2.5275,0.14)(2.7,0.085)
\pscircle(2.5,0.15){0.03}
\psdot(0.67,0.7)
\put(0.67,0.78){$p$}
\psdot(1.29,0.515)
\put(1.28,0.39) {$a$}
\psdot(0.9,0.15)
\put(0.88,0.04) {$q$}
\psline[linestyle=dashed](0.67,0.7)(0.9,0.15)
\psdot(1.49,0.7)
\put(1.5,0.78) {$r$}
\qline(0.9,0.15)(1.49,0.7)
\put(-0.15,0.12) {$L$}
\put(-0.15,0.68) {$K$}
\put(-0.15,0.88) {$M$}
\endpspicture}

\def\EuclidParallelRawFigure{%
\pspicture(3.2, 1.2)
\pspolygon[fillstyle=solid,fillcolor=lightgray](0.67,0.7)(0.78,0.425)(1.49,0.7)
\pspolygon[fillstyle=solid,fillcolor=lightgray](0.07,0.15)(0.78,0.425)(0.9,0.15)
\qline(-0.15,0.15)(2.47,0.15)
\qline(2.53,0.15)(2.7,0.15)
\qline(-0.15,0.7)(2.7,0.7)
\qline(0.0,0.9)(2.473,0.163)
\qline(2.5275,0.14)(2.7,0.085)
\pscircle(2.5,0.15){0.03}
\put(2.5,0.04){$x$}
\psdot(0.67,0.7)
\put(0.67,0.78){$p$}
\psdot(1.29,0.515)
\put(1.28,0.39) {$a$}
\psdot(0.9,0.15)
\put(0.88,0.04) {$q$}
\psdot(0.07,0.15)
\put(0.07,0.04){$s$}
\qline(0.67,0.7)(0.9,0.15)
\psdot(1.49,0.7)
\put(1.5,0.78) {$r$}
\qline(1.49,0.7)(0.07,0.15)  
\psdot(0.78,0.425)
\put(0.71,0.3){$t$}
\qline(0.9,0.15)(1.49,0.7)  
\put(-0.3,0.12) {$L$}
\put(-0.3,0.68) {$K$}
\put(-0.15,0.88) {$M$}
\endpspicture}

\title{Herbrand's theorem and non-Euclidean geometry}
\date{September 17, 2012}
 
\author{Michael Beeson}
\address{San Jos\'e State University}
\author{Pierre Boutry}
\address{ICube, UMR 7357 CNRS, University of Strasbourg}
\author{Julien Narboux}
\address{ICube, UMR 7357 CNRS, University of Strasbourg}

\begin{document}
\pagestyle{headings}  

\maketitle

\begin{abstract}  
We use Herbrand's theorem to give a new proof that Euclid's 
parallel axiom is not derivable from the other axioms of first-order
Euclidean geometry.  Previous proofs involve constructing models of  
non-Euclidean geometry.  This proof uses a very old and basic 
theorem of logic together with some simple properties of ruler-and-compass
constructions to give a short, simple, and intuitively appealing proof.
\end{abstract}

\section{Introduction}  We intend this paper to be read by mathematicians
who are unfamiliar with mathematical logic and also unfamiliar with 
non-Euclidean geometry;  therefore we ask the patience of readers who 
are familiar with one or both of these subjects.  

 We begin with a brief discussion of axioms 
for plane Euclidean geometry.  Every such axiom system will have variables
for points.   Some axiom systems may have variables
for other objects, such as lines or angles, but Tarski showed that these
are not really necessary.  For example, angles can be discussed in terms of 
ordered triples of points, and lines in terms of ordered pairs of points.
For simplicity we focus on such a points-only axiomatization.

The primitive relations of such a theory usually include a ``betweenness''
relation, and an ``equidistance'' relation.   We write
$\T(a,b,c)$ to express that $b$ lies (non-strictly)
between $a$ and $c$ (on the same line), and $\E(a,b,c,d)$ to express that segment
$ab$ is congruent to segment $cd$.  $\E$ stands for ``equidistance'', because in the 
standard model ``congruent'' means that the distance $ab$ is equal to the distance $cd$;
but there is nothing in the axioms about numbers to measure distance, or about distance
itself.  Sometimes it is convenient
to use $\B(a,b,c)$ for strict betweenness, i.e. $a \neq b$ and $b \neq c$ and $\T(a,b,c)$.

Some of the axioms will assert the existence of ``new'' points that are 
constructed from other ``given'' points in various ways.  For example, 
one axiom says that segment $ab$ can be extended past $b$ to a point $x$,
lying on the line determined by $ab$, such that segment $bx$ is congruent to 
a given segment $pq$.   That axiom can be written formally, using the logician's
symbol $\land$ for ``and'', as 
$$ \exists x\,( \T(a,b,x) \land \E(b,x,p,q))$$
It is   possible to replace the quantifier $\exists$ with a ``function symbol''. 
We denote the point $x$ that is asserted to exist by $ext(a,b,p,q)$.  Then the 
axiom looks like
$$  \T(a,b,ext(a,b,p,q)) \land \E(b,ext(a,b,p,q),p,q)$$
This transformation is called \textit{Skolemization}.
This form is called ``quantifier-free'',  because $\exists$ and $\forall$ are 
called ``quantifiers'',  and we have eliminated the quantifiers. Although the 
meaning of the axioms is the same as if it had $\forall a,b,p,q$  in front,
the $\exists$ has been replaced by a function symbol.

When a theory has function symbols, then they can be combined. For example,
$ext(a,b,ext(u,v,p,q),ext(a,b,p,q))$  is a term.  The definition of ``term''
is given inductively:  variables are terms,
constants are terms, and if one substitutes
terms in the argument places of function symbols, one gets another term.  

In Tarski's axiomatization of geometry,  there are only a few axioms that 
are not already quantifier-free.  One of them is the segment extension axiom 
already discussed.  Another is Pasch's axiom.   Moritz Pasch originally proposed
this axiom in 1852, to repair the defects of Euclid.   It intuitively says that
if a line meets one side of a triangle and does not pass through the endpoints of 
that side, then it must meet one of the other sides of the triangle.  In other 
words, under certain circumstances, there will exist the intersection point of two lines.
A quantifier-free version of Tarski's axioms will contain a function symbol for 
the point asserted to exist by (a version of) Pasch's axiom.  

Another axiom in Tarski's theory asserts the existence of an intersection point of  
a circle and a line, provided the line has a point inside and a point outside the circle.
Another function symbol can be introduced for that point.  Then the terms of this 
theory correspond to certain ruler-and-compass constructions.  The number of symbols
in such a term corresponds to the number of ``steps'' required with ruler and compass 
to construct the point defined by the term.

The starting point for the work reported here is this:  a quantifier-free theory 
of geometry, whose terms correspond to ruler-and-compass constructions, viewed 
as a special case of situation of much greater generality:  some first-order, 
quantifier-free theory.   Herbrand's theorem applies in this much greater generality,
and we will simply investigate what it says when specialized to geometry.

\section{Herbrand's theorem}
Herbrand's theorem is a general logical theorem about any axiom system whatsoever
that is 
\begin{itemize}
\item first-order, i.e. has variables for some kind(s) of objects, but not 
for sets of those objects,  and 
\item quantifier-free, i.e. $\exists$ has been replaced by function symbols
\end{itemize}
Herbrand's theorem says that under these assumptions, if the theory proves
an existential theorem $\exists y\, \phi(a,y)$,  with $\phi$ quantifier-free,
 then there exist finitely many
terms $t_1,\ldots,t_n$ such that the theory proves
$$  \phi(a,t_1(a)) \lor \phi(a,t_2(a)) \ldots \lor \ldots \phi(a,t_n(a)).$$
The formula $\phi$ can, of course, have more variables that are not explicitly shown here,
and $a$ and $x$ can each be several variables instead of just one, in which case the $t_i$ stand for 
corresponding lists of terms.  For a proof see \cite{buss-handbook}, p.~48. 

In order to illustrate the theorem, consider the example when $\phi$ is $\phi(a,b,c,x,y)$,
and it says that $a \neq b$, and $x$ lies on the line determined by $ab$, and $y$ does 
not lie on that line, and $xy$ is perpendicular to $ab$ and $c$ is between $x$ and $y$.
   Collinearity can be expressed
using betweenness, and the relation $xy \perp ab$ can also be expressed using betweenness
and equidistance.   Then $\exists x,y\,\phi(x,y)$ says that there exists a line through 
point $c$ perpendicular to $ab$.  Usually in geometry, we give two different constructions 
for such a line, according as $c$ lies on line $ab$ or not.  If it does, we ``erect'' a 
perpendicular at $c$, and if it does not, we ``drop'' a perpendicular from $c$ to line $ab$.
\def\foot{{\mbox{\it foot}}}
When we ``drop'' a perpendicular, we compute $\foot_1(a,b,c)$, and we can define $head_1(a,b,c) = c$.
When we ``erect'' a perpendicular, we compute $head_2(a,b,c)$, and we can define
$\foot_2(a,b,c)=c$.
Thus if $c$ is not on the line, we have $\phi(a,b,c,\foot_1(a,b,c),head_1(a,b,c))$,
and if $c$ is on the line, we have $\phi(a,b,c,\foot_2(a,b,c),head_2(a,b,c)$.  Since 
$c$ either is or is not on the line we have
$$ \phi(a,b,c,\foot_1(a,b,c),head_1(a,b,c)) \lor \phi(a,b,c,\foot_2(a,b,c),head_2(a,b,c)) $$
Comparing this to Herbrand's theorem, we see that we have specifically constructed 
examples of two lists (of two terms each)   $t_1$ and $t_2$  illustrating that Herbrand's
theorem holds in this case.   Herbrand's theorem, however,  tells us {\em without doing any 
geometry}  that if there is any proof at all of the existence of a perpendicular to $ab$ through $c$, 
from the axioms of geometry mentioned above, 
then there must be a finite number of ruler-and-compass constructions such that, for every 
given $a,b,c$, one of those constructions works.   We have verified, using geometry,
that we can take the ``finite number'' of constructions to be 2 in this case,  but the 
beauty of Herbrand's theorem lies in its generality.

\section{Non-Euclidean geometry}
\begin{figure}[ht]
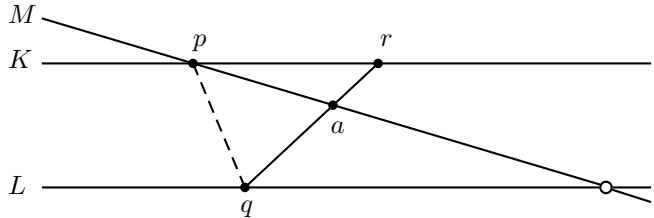
   
\caption{Euclid 5.  $M$ and $L$ must meet on the right side, provided $\B(q,a,r)$ and $pq$ makes
alternate interior angles equal with $K$ and $L$.
\label{figure:EuclidParallelFigure}
}
\hskip 2.5cm
\EuclidParallelFigure
\bigskip
\end{figure}

Euclid listed five axioms or postulates, from which, along with his ``common notions'',
he intended to derive all his theorems.  The fifth postulate, known as ``Euclid 5'',
had to do with parallel lines, and is also known as the ``parallel postulate.''
See Fig.~\ref{figure:EuclidParallelFigure}.

From antiquity, mathematicians felt that Euclid 5 was less ``obviously true'' than the other 
axioms, and they attempted to derive it from the other axioms.  Many false ``proofs'' were 
discovered and published.  All this time, mathematicians felt that geometry was ``about''
some true notion of space,  which was either given by the physical space in which we live,
or perhaps by the nature of the human mind itself.  Finally, after constructing long chains 
of reasoning from the assumption that the parallel postulate is false, some people came to 
the realization that there could be ``models of the axioms'' in which ``lines'' are 
interpreted as certain curves, and ``distances'' also have an unusual interpretation.  Such 
models were constructed in which Euclid 5 is false, but the other axioms are true.   Hence,
Euclid 5 can never be proved from the other axioms.  There was a good reason for all those 
failures!   See \cite{greenberg} and \cite{hartshorne} for the full history of these fascinating 
developments, and descriptions of the models in question.

\section{Tarski's axioms for geometry} 
In order to state our theorem precisely, we need to mention a specific axiomatization of 
geometry.  For the sake of definiteness, we use the axioms (A1-A11) of Tarski, as set forth 
in the definitive reference \cite{schwabhauser}. We list those axioms in Table~\ref{tarski-axioms}. Those who do not read German can consult 
 \cite{tarski-givant}.
 
\begin{table}
\caption{Tarski's axioms for geometry\label{tarski-axioms}}
\begin{tabular}{cll}
A1 & Symmetry & $\congT{A,B}{B,A}$\\
A2 & Pseudo-Transitivity\qquad \qquad & $\congT{A,B}{C,D} \land \congT{A,B}{E,F} \imp \congT{C,D}{E,F}$\\
A3 & Cong Identity & $\congT{A,B}{C,C} \imp a=b$\\
A4 & Segment extension &$\exists e( \bT{A}{B}{E} \land \congT{B,E}{C,D})$\\
A5 & Five segments&$\begin{array}[t]{l}
 \congT{A,B}{A',B'} \land \congT{B,C}{B',C'} \land \\ \congT{A,D}{A',D'}
  \land \congT{B,D}{B',D'} \land a \neq b  \land\\
\bT{A}{B}{C} \land \bT{A'}{B'}{C'} \imp \congT{C,D}{C',D'}\\
\end{array}$\\
A6 & Between Identity & $\bT{A}{B}{A} \imp a=b$\\
A7 & Inner Pasch & $\begin{array}[t]{l}
             \bT{A}{P}{C} \land \bT{B}{Q}{C} \imp \\
             \exists x\,(\bT{P}{X}{B} \land \bT{Q}{X}{A})\end{array}$\\
A8 & Lower Dimension &$\exists a b c( \lnot \bT{A}{B}{C} \land \lnot \bT{B}{C}{A} \land \lnot \bT{C}{A}{B} )$\\
A9 & Upper Dimension &$ \congT{A,P}{A,Q} \land \congT{B,P}{B,Q} \land \congT{C,P}{C,Q} \land p \neq q  $\\
& & $ \imp\bT{A}{B}{C} \lor \bT{B}{C}{A} \lor \bT{C}{A}{B}$\\
A10 & Parallel & $\exists x y( \bT{A}{D}{T} \land \bT{B}{D}{C} \land a \neq d \imp$\\ & & $ \bT{A}{B}{X} \land \bT{A}{C}{Y} \land \bT{X}{T}{Y}$)\\
A11 & Continuity & $\forall X Y
    ((\exists a (\forall x y, x \in X \land y \in Y \imp \bT{A}{x}{y}))$\\
  &  &$ \imp\exists b (\forall x y, x \in X \land y \in Y \imp \bT{x}{B}{y}))$\\
CA & Circle axiom & $\begin{array}[t]{l} \bT{a}{x}{b} \land \bT{a}{b}{y} \land \congT{a,x}{a,p} \land \\ 
\congT{a,q}{a,y} \imp \exists z( \congT{a,z}{a,b} \land \bT{p}{z}{q}) \end{array}$
\end{tabular}
\end{table}

 Of these axioms, we need concern ourselves in detail only with those few that are not already quantifier-free.  Axiom (A4) is the segment extension axiom discussed above; we introduce the symbol $ext(a,b,p,q)$ to express it in quantifier-free form. 
 The lower dimension axiom (A8) states that there exists three non collinear points. We introduce three constants
$\alpha$, $\beta$, and $\gamma$   to express it in quantifier-free form.  The two 
modified axioms are explicitly:
\begin{table}
\caption{Axioms A4 and A8 in quantifier-free form}
\begin{tabular}{cll}
A4${}^{\prime}$ & Segment extension \qquad\qquad&$   \bT{A}{B}{E} \land \congT{B,E}{C,D})$\\
A8${}^{\prime}$ & Lower Dimension &$  \lnot \bT \alpha \beta \gamma \land \lnot \bT \beta \gamma \alpha \land \lnot \bT \gamma\alpha\beta )$
\end{tabular}
\end{table}

\subsection{Pasch's axiom}  
Moritz Pasch \cite{pasch1882} (see also \cite{pasch1926}, with an historical appendix by Max Dehn)
 supplied (in 1882) an axiom that repaired many of the defects that 
nineteenth-century rigor found in Euclid.  Roughly, a line that enters a 
triangle must exit that triangle.
As Pasch formulated it,  it is not in $\forall\exists$ form.  There are two $\forall\exists$  versions,
illustrated in Fig.~\ref{figure:pasch}.  These formulations of Pasch's axiom 
go back to Veblen \cite{veblen1904}, who proved outer Pasch implies inner Pasch. 
Tarski originally took outer Pasch as an axiom.  In \cite{gupta1965}, Gupta proved
both that inner Pasch implies outer Pasch, and that outer Pasch implies inner Pasch,
using the other axioms of the 1959 system.  In the final version \cite{schwabhauser},
inner Pasch is an axiom.  Here are the precise statements of the axioms illustrated in Fig.~\ref{figure:pasch}:

\begin{figure}[ht]
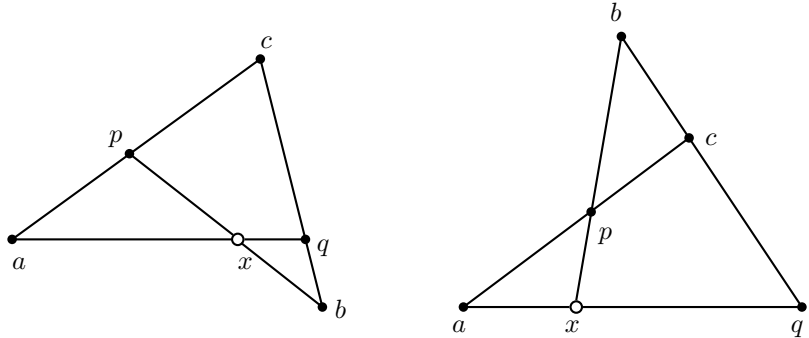
  
\caption{ Inner Pasch (left) and Outer Pasch (right).  Line $pb$ meets triangle $acq$ in one side.
 The open circles show the points asserted to exist on the other side.}
\label{figure:pasch}
\vskip -0.2cm
\hskip 0.7cm
\InnerOuterPaschFigure
\vskip -0.8cm
\end{figure}
\smallskip

\axioms
$\T(a,p,c) \land \T(b,q,c) \implies \exists x\,(\T(p,x,b) \land \T(q,x,a))$\qquad\qquad& (A7) inner Pasch \\
$\T(a,p,c)\land \T(q,c,b) \implies \exists x\,(\T(a,x,q) \land \T(b,p,x))$&  outer Pasch 
\endaxioms
\smallskip

\noindent
In order to express inner Pasch in quantifier free form, we introduce the symbol $ip(a,p,c,b,q)$ 
for the point $x$ asserted to exist.   This corresponds to the ruler-and-compass
(actually just ruler) construction of finding the intersection point of lines $aq$ and $pb$.
There is a codicil to that remark, in that Tarski's axiom allows the degenerate case in which 
the segments $aq$ and $pb$ both lie on one line (so that there are many intersection points,
rather than a unique one), but we do not care in this paper that in such a case the 
construction cannot really be carried out with ruler and compass.  
Also, we call the reader's attention to this fact:  point $c$ is not needed to draw the lines with a ruler, but it is needed to ``witness'' 
that the lines actually ``should'' intersect.

\subsection{Continuity and the Circle Axiom}
\begin{figure}[ht]
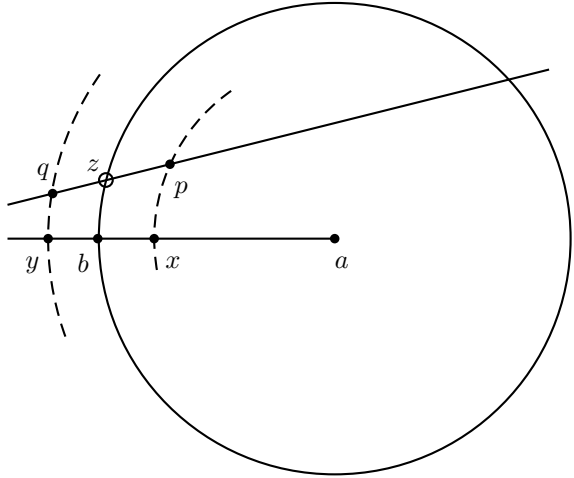

\caption{Circle Axiom(CA). Point $p$ is inside, $q$ is outside, so $pq$ meets the circle.
\label{ca}}
\begin{center}
\SegmentCircleContinuityFigure
\bigskip
\end{center}
\end{figure}

Axiom (A11) is the ``continuity'' axiom.  In its full generality, it says that 
``first-order Dedekind cuts are filled.''  Closely related to (A11) is the 
``circle axiom'' (CA),  
which says that if $p$ lies inside the circle with center $a$ and passing through $b$, 
and $q$ lies outside that circle, then segment $pq$ meets the circle (see Fig.~\ref{ca}).%
\footnote{There is no ``standard'' name for this axiom. 
Tarski did not give the this axiom a name, only a number; in \cite{schwabhauser} and other 
German works it is called the ``Kreisaxiom'', which we translate literally here.  In
\cite{gupta1965} it is called the ``line and circle intersection axiom'', which we find too long.
In \cite{greenberg} (p.~131)  it is called the ``segment-circle continuity principle.'' 
}

Points $x$ and $y$ in the figure serve as ``witnesses'' that $p$ and $q$ are inside and outside,
respectively.
Specifically,  ``$p$ lies inside
the circle''  means that $ap < ab$,  which in turn means that there is a point $x$ between $a$ and $b$
such that $\E(a,x,a,p)$, i.e. segment $ax$ is congruent to $ap$.  Similarly, ``$q$ lies outside
the circle'' means  there exists $y$ with $\B(a,b,y)$ and $\E(a,q,a,y)$.  In order to express 
segment-circle continuity in quantifier-free form, we can introduce
a symbol $\ilc(p,q,a,b,x,y)$ for the point of intersection of $pq$ with the circle.   Even though 
$x$ and $y$ are not needed for the ruler-and-compass construction of this point, they must be 
included as parameters of $\ilc$.  

We return below to the general axiom (A11), but first we show how to finish the proof of our 
main theorem if only 
the circle axiom  is used, instead of the full schema (A11).

\subsection{The parallel axiom}
\begin{figure}[h]
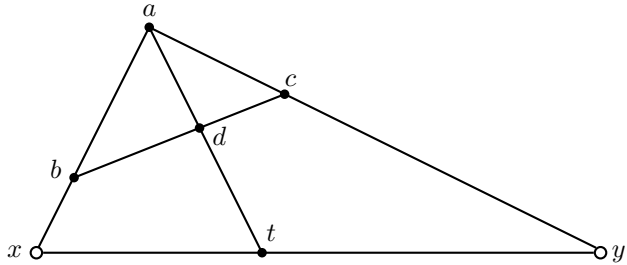
   
\caption{Tarski's parallel axiom (A10).}
\label{figure:TarskiParallelFigure} 
\hskip 1in 
\TarskiParallelFigure
\end{figure}

Tarski used a variant formulation (A10) of Euclid 5, illustrated in Fig.~\ref{figure:TarskiParallelFigure}.
One can prove the equivalence (A10) with Euclid 5, and (A10) has the advantage of being 
very simply expressed in a points-only language.  Open circles indicate the two points asserted
to exist.

\begin{figure}[ht]
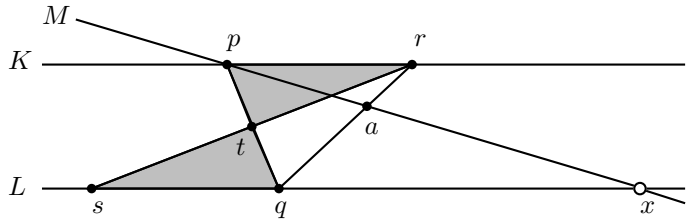

\caption{Euclid 5.  Transversal $pq$ of lines $M$ and $L$ makes corresponding interior angles less than 
two right angles, as witnessed by $a$. The shaded triangles are assumed congruent. Then $M$ 
meets $L$ as indicated by the open circle.}
\label{figure:EuclidParallelRawFigure}
\hskip 2.5cm
\EuclidParallelRawFigure
\end{figure}

For our independence proof, we work with Tarski's axiom A10 rather than with Euclid 5.
Nevertheless, we include a
formulation of Euclid's parallel postulate, expressed in 
Tarski's  language.  Euclid's version mentions angles, and the 
concept of ``corresponding interior angles'' made by a transversal.  Fig.~\ref{figure:EuclidParallelRawFigure} illustrates
the following points-only version of Euclid 5.    

\medskip

\axioms
$   \B(q,a,r) \land \B(p,t,q) \land pr=qs \land pt=qt \land rt=st      $&(Euclid 5)\\
$\land\  \neg\, Col(s,q,p)  \implies \exists x\,( \B(p,a,x) \land \B(s,q,x))$&
\endaxioms
\smallskip

\section{Consistency of non-Euclidean geometry via Herbrand's theorem}
The point of this paper is to show that one can use the very general theorem of Herbrand
to prove the consistency of non-Euclidean geometry, doing extremely little actual geometry.
All the geometry required is the observation that 
when we construct points from some given points, at each construction 
stage the maximum distance between the points at most doubles.  

In order to state our theorem precisely, we define 
 $T$ to be Tarski's ``neutral ruler-and-compass geometry'',  where ``neutral''
means that the parallel axiom (A10) (equivalent to Euclid 5)  
is not included, and ``ruler-and-compass'' means that (A11)
is replaced by the circle axiom (CA).  In addition, $T$ uses the 
quantifier-free versions of the segment-extension and dimension axioms discussed above.
The following 
lemma states precisely what we mean by, ``at each construction state the maximum distance
between the points at most doubles.''

\begin{Lemma} \label{lemma:distance}
 The function symbols of $T$ have the following property, when 
interpreted in the Euclidean plane $\R^2$:   if $f$ is one of those function symbols,
i.e. $f$ is $ext$ or $\ilc$ or $ip$,  then the distance of $f(x_1,\ldots,x_j)$ from 
any of the parameters $x_1,\ldots x_j$ is bounded by twice the maximum distance between 
the $x_j$. 
\end{Lemma}

\noindent{\em Proof}.
When we extend a segment $ab$ by a distance $pq$,  the distance of the new point $ext(a,b,p,q)$
from the points $a,b,p,q$ is at most twice the maximum of $ab$ and $pq$.  
The point constructed by $ip$ 
is between some already-constructed points, so $ip$ does not increase the distance at all.
The point constructed by $\ilc$ is no farther
from the center $a$ of the circle than the given point $b$ on the circle is, and hence no more than 
$ab$ farther from any of the other points, and hence no more than twice as far from any of the
other parameters of $\ilc$ 
 as the maximum distance between those points.

\begin{Theorem}  Let $T$ be Tarski's ``neutral ruler-and-compass geometry'',  where ``neutral''
means that the parallel axiom (A10) (equivalent to Euclid 5)  
is not included, and ``ruler-and-compass'' means that (A11)
is replaced by the ``circle axiom'' (CA).   Then $T$ does not prove the parallel axiom (A10).
\end{Theorem}

\begin{figure}
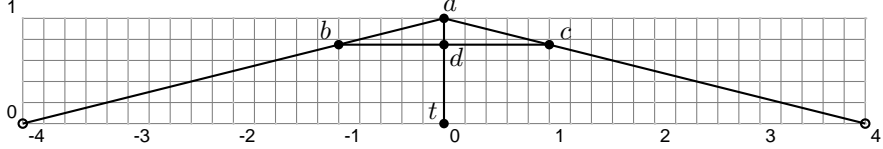

\caption{Construction of a point too far away. Here $k=2$ and the constructed  points are indicated by the open circles.}
\label{figure:toofar}
\begin{center}
\CounterExampleFigure
\end{center}
\end{figure}

\noindent{\em Proof}.  Suppose, for proof by contradiction, that $T$ does prove (A10).
There is a formula $\phi(a,b,c,d,t,x,y)$ such that axiom (A10) has the form
$$ \exists x,y\, \phi(a,b,c,d,t,x,y),$$
where $\phi$ expresses the betweenness relations shown in the figure.
Then, by Herbrand's theorem, there are finitely many 
terms $X_i(a,b,c,d,t)$ and $Y_i(a,b,c,d,t)$, for $i=1,2,\ldots,n$, such that 
$T$ proves
$$ \bigvee_{i=1}^n \phi(a,b,c,d,t,X_i(a,b,c,d,t), Y_i(a,b,c,d,t)).$$
Let $k$ be an integer greater than the maximum number of function symbols in any of those $2n$ terms.
Choose points $a,b,c,d$ and $t$ in the ordinary plane $\R^2$ as follows (see Fig.~\ref{figure:toofar})
\begin{eqnarray*}
t &=& (0,0) \\
a &=& (0,1) \\
b &=& (-1,1-2^{-k-2}) \\
c &=& (1,1-2^{-k-2}) \\
d &=& (0,1-2^{-k-2}) 
\end{eqnarray*}

Suppose $x$ and $y$ are as in (A10); then one of them has a nonnegative second coordinate,
and the other one must have a first coordinate of magnitude at least $2^{k+2}$.   But 
then, according to the lemma, it cannot be the value of one of the terms $X_i(a,b,c,d,t)$
or $Y_i(a,b,c,d,t)$,  which, since they involve $k$ symbols starting with points no more 
than distance 2 apart, cannot be more than $2^{k+1}$ from any of the starting points.
This contradiction completes the proof.

\section{Full first-order continuity}
In this section we show how to extend the above proof to include the full (first-order)
continuity 
axiom (A11) instead of just the circle axiom.  The difficulty is that (A11)
is far from quantifier-free, but instead is an axiom schema.  That means, it is actually 
an infinite number of axioms,  one for each pair of first-order formulas $(\phi,\psi)$.
The axiom says, if the points satisfying $\phi$ all lie on a line to the left of the points
satisfying $\psi$, then there exists a point $b$ non-strictly between any pair of points
$(x,y)$ such that $\phi(x)$ and $\psi(y)$.  

The keys to extending our proof are Tarski's deep theorem on quantifier-elimination for 
algebra,  and the work of Descartes and Hilbert on defining arithmetic in geometry.  
Modulo these results, which in themselves have nothing to do with non-Euclidean geometry,
the proof extends easily to cover full continuity, as we shall see.

A {\em real-closed field} is an ordered field $\F$ in which every polynomial of odd degree
has a root, and every positive element has a square root.  Tarski proved in \cite{tarski1951} the following 
fundamental facts:

\begin{itemize} 
\item  Every formula in Tarski's language is provably equivalent to a 
quantifier-free formula. 
\item Every model of Tarski's axioms has the form $\F^2$, where $\F$ is a real-closed
field, and betweenness and equidistance are interpreted as you would expect.
\end{itemize}

Since Descartes and Hilbert showed how to give geometric definitions of addition, multiplication,
and square root, there are formulas in Tarski's language defining the operations of multiplying
and adding points on a fixed line $L$, with points $0$ and $1$ arbitrarily chosen on $L$,
 and taking square roots of points to the right of $0$ (see chapter 14 and 15 of \cite{schwabhauser}).   Since the existence of square roots
 follows from the circle axiom,  the full continuity schema is equivalent 
to the schema that expresses that polynomials of odd degree have zeroes:

\begin{equation}
  \exists x\, (a_0 + a_1 x + \ldots + a_{n-1} x^{n-1} + x^{n} = 0 ). \label{eq:1}
\end{equation}
Note that without loss of generality the leading coefficient can be taken to be 1.
Here the algebraic notation is an abbreviation for geometric formulas in Tarski's language.
The displayed formula represents one geometric formula for each fixed odd integer $n$, so it 
still represents an infinite number of axioms, but Herbrand's theorem applies even if there 
are an infinite number of axioms.   The essential point is that this axiom schemata is 
purely existential, so we can make it quantifier-free by introducing a single new function 
symbol $f(a_0,\ldots,a_{n-1})$ for a root of the polynomial.

\begin{Theorem}  Axioms A1-A9 and axiom schema A11 together do not prove the 
parallel axiom A10.
\end{Theorem}

\noindent{\em Proof}.  Suppose, for proof by contradiction, that A10 is provable from A1-A9 
and A11.   Then, the models of A1-A9 and A11 are all isomorphic to planes over real-closed fields.
Then, as explained above, the full schema A11 is equivalent (in the presence of A1-A10) 
to the schema (\ref{eq:1}) plus
the circle axiom.%
\footnote{It is worth emphasizing that this equivalence depends on developing the theory 
of perpendiculars without any continuity axiom at all, not even the circle axiom.  This 
was one of the main results of \cite{gupta1965}, and is presented in \cite{schwabhauser}, where
it serves as the foundation to the
 development of arithmetic in geometry.  It is quite difficult even to prove the 
circle axiom directly from A11 without Gupta's results,  although Tarski clearly believed
decades earlier that the circle axiom does follow from A1--A11, or he would have included 
it as an axiom.}
  
That is,  it suffices to 
 supplement ruler-and-compass constructions by the ability 
to take a root of an arbitrary polynomial.  The point that allows our proof to work is 
simply that the roots of polynomials can be bounded in terms of their coefficients.  For example,
the well-known ``Cauchy bound'' says that any root is bounded by the maximum of $1 + \vert a_i \vert$ 
for $i=0,1,\ldots n-1$, which is at most 1 more than the max of the parameters of $f(a_0,\ldots,a_{n-1})$.
Below we give, for completeness, a short proof of the Cauchy bound, but first, we finish the proof of the theorem.

We can then modify Lemma~\ref{lemma:distance} to say that the distance is at most the max of 1 and double the previous 
distance.  In the application we start with points that are 1 apart, so the previous argument
applies without change.  That completes the proof.

\begin{Lemma}[Cauchy bound] The real roots of $a_0 + a_1x + \ldots + a_{n-1}x^{n-1} + x^n$
are bounded by the maximum of $1+ \vert a_i\vert $.
\end{Lemma}

\noindent{\em Proof}.   Suppose $x$ is a root.  If $\vert x \vert \le 1$ then $x$ is bounded, hence we may assume 
$\vert x \vert > 1$.    Let $h$ be the max of the $|a_i|$.  Then 
\begin{eqnarray*}
-x^n &=& \sum_{i=0}^{n-1} a_ix^i , \ \ \mbox{so \ \ } 
\vert x \vert^n \le  h \sum_{i=0}^{n-1} \vert x \vert^i \ = \  h \frac { \vert x \vert^n -1 } {\vert x \vert -1}
\end{eqnarray*}
Since $\vert x \vert > 1$ we have
\begin{eqnarray*}
\vert x \vert -1 \le h\frac { \vert x \vert^n -1 } {\vert x \vert^n} \ \le h. 
\end{eqnarray*}
Therefore $\vert x \vert \le 1 + h$.
That completes the proof.

\section{Related proof-theoretical work of others}
Skolem \cite{skolem1920} already in 1920 proved the independence of a form of the parallel axiom 
from the other axioms of projective geometry, using methods similar to Herbrand's theorem.
In 1944, Ketonen invented the system of sequent calculus made famous in Kleene \cite{kleene1952}
as G3, and used it to reprove Skolem's result and extend it to affine geometry.  This result
was reproved using a different sequent calculus in 2001 by von Plato \cite{vonplato2001}. 
It should be noted that the modern proof of Herbrand's theorem also proceeds by cut-elimination
in sequent calculus.  Our proof of the independence of Euclid's parallel axiom improves on 
these past results in that (i)  it works for ordinary geometry, not just for projective or affine geometry, and (ii) it depends on proof theory only for Herbrand's theorem:  no direct analysis
or even mention of cut-free proofs is required.

\section{Another proof via a model of Max Dehn's}
Max Dehn, a student of Hilbert, gave a model of A1-A9 plus the circle axiom. 
 Dehn's model is easily described and, like our proof,
has no direct relationship to non-Euclidean geometry.  

An element $x$ in an ordered field $\K$ is called {\em finitely bounded} if it is less than 
some integer $n$, where we identify $n$ with $\sum_{k=1}^n 1$.   $\K$ is {\em Archimedean} if 
every element is finitely bounded.  It is a simple exercise to construct a non-Archimedean
Euclidean field, or even a non-Archimedean real-closed field. (For details 
about Dehn's model, see  Example 18.4.3  and Exercise 18.4
of \cite{hartshorne}.)
Dehn's model begins with a non-Archimedean Euclidean field $\K$.  Then the set $\F$ of finitely 
bounded elements of $\K$ is   a Euclidean ring, but not a Euclidean field: there are elements
$t$ such that $1/t$ is not finitely bounded.   These are called ``infinitesimals.''  
Dehn's point was that $\F^2$ still satisfies the axioms of ``Hilbert planes'', which are 
equivalent (after \cite{schwabhauser}) to A1-A9.  The reason is similar to the reason that 
our Herbrand's-theorem proof works:  the constructions given by segment extension and 
Pasch's axiom can at most double the size of the configuration of constructed points, 
so they lead from finitely bounded points to other finitely bounded points.  
  Since square roots of finitely bounded 
elements are also finitely bounded, $\F^2$ satisfies the circle axiom too.  But
$\F^2$ does not satisfy the parallel axiom, since there are lines with infinitesimal slope
through $(0,1)$ that do not meet the $x$-axis of $\F$.  (They meet the $x$-axis of $K$,
but not at a finitely bounded point.)   
  
In this way Dehn showed that (the Hilbert-style equivalent of) A1-A9, together with 
the circle axiom, does not imply the parallel postulate A10.  We add to Dehn's
  proof the extension to the full first-order continuity schema A11,  by the same 
  trick as we used for our Herbrand's-theorem proof.  Namely, suppose for proof 
  by contradiction that A10 is provable
  from A1-A9 and A11.  Then in A1-A9 plus segment-circle continuity, A11 is equivalent 
  to the schema (\ref{eq:1}) saying that odd-degree polynomials have roots.   Now  
  construct Dehn's model starting from a non-Archimedean real-closed field $\K$. 
Then $\F$ still satisfies (\ref{eq:1}), because of the Cauchy bound:  if the coefficients
$a_i$ are finitely bounded, so are the roots of the polynomial.   But then $\F^2$
satisfies A11,  and hence, according to our assumption, it satisfies A10 as well; but 
we have seen that it does not satisfy A10, so we have reached a contradiction.  That 
contradiction shows that A10 is not provable from A1-A9 and A11.

Note that this proof, like the proof via Herbrand's theorem,  does not actually construct a model 
of non-Euclidean geometry, that is, a model satisfying A1-A9, A11, but not A10.  
That is the interest of both proofs:  the consistency of non-Euclidean geometry is shown,
in the one case by proof theory, and the other by algebra (or model theory if you prefer to 
call it that), without doing any non-Euclidean geometry at all.   Moreover, the classical 
constructions of models of non-Euclidean geometry (the Beltrami-Klein and Poincar\'e models 
described in \cite{greenberg}, Ch.~7), satisfy not only the first-order continuity schema 
but also the full second-order continuity axioms.  Herbrand's theorem is about first-order
logic, so it cannot replace these classical geometrical constructions; but still, we
have shown here that a little logic goes a long ways.

\nocite{cavinessjohnson}
\bibliographystyle{asl}
\bibliography{HerbrandEuclid}
\end{document}